\documentclass{amsart}
\usepackage{amsmath}
\usepackage{mathtools}

\usepackage{pdfpages}

\theoremstyle{theorem}

\usepackage{tikz}
\usepackage{pgfplots}
\usepackage{hyperref}

\usepackage[bibtex-style]{amsrefs}
\usepackage{amsmath, amsfonts, amssymb, amsthm, amscd}

\begin{document}

\title{Using binary string to prove the Collatz conjecture}
\markright{Prove the Collatz conjecture}
\author{Ji she Feng}

%\small {}

\footnotetext{School of Mathematics and Information Engineering, Longdong University,  Qingyang,  Gansu,  745000,
 China E-mail: gsqyfjsfjs6567@sina.com.}
\footnotetext{Foundation item: Educational technology innovation project of Gansu Province (No. 2022A-133)}

\begin{abstract}
We introduce a full binary directed tree structure to represent the set of natural numbers, further categorizing them into three distinct subsets: pure odd numbers, pure even numbers, and mixed numbers. We adopt a binary string representation for natural numbers and elaborate on the composite methodology encompassing odd- and even-number functions. Our analysis focuses on examining the iteration sequence (or composition) of the Collatz function and its reduced variant, which serves as an analog to the inverse function, to scrutinize the validity of the Collatz conjecture. To substantiate this conjecture, we incorporate binary strings into an algebraic formula that captures the essence of the Collatz sequence. By this means, we transform discrete powers of 2 into continuous counterparts, ultimately culminating in the smallest natural number, 1. Consequently, the sequence generated through infinite iterations of the Collatz function emerges as an eventually periodic sequence, thereby validating an enduring 87-year-old conjecture.
\end{abstract}

\maketitle

\section{Introduction} % (fold)
\label{sec:introduction}

In the study of number theory, odd and even numbers are a fundamental pair of ideas. Natural number sets can be parted to two different sets. There are many conjectures that attempt to generalize the fact of different kinds of natural numbers discovered in a restricted range to the entire infinite set of natural numbers. This article will proof the famous Collatz conjecture, which states that for each natural number $n$, if it is even, divided by 2, if it is odd, multiplied by 3, and added 1, and so on, the eventual value must be 1. It is also referred to as the $3n+1$ conjecture and was put forth in 1937 by Lothar Collatz, also known as the $3n+1$ problem. The mathematician Paul Erdos once said of this conjecture: "Mathematics may not be ready for such problems"$^{[1,2]}$.

Inspired by Euler's dot-line graph in graph theory for solving the Konigsberg Seven Bridge problem, we have confidence that a similar solution can be found. Leveraging our knowledge of piecewise and iteration functions, binary strings, and full binary directed trees, we utilize binary strings to illustrate the step-by-step progression of odd-function and/or even-function iterations, concealed within the Decimal Number System. While no new mathematical concepts are introduced, we are convinced that the existing framework is adequate to substantiate the conjecture.

For the Collatz conjecture, we can describe it as a piecewise function:
\begin{equation}
T(n)=\left\{
\begin{array}{cc}
3n+1, & \text{if }n\text{ is odd number,} \\
\frac{n}{2} & \text{if }n\text{ is even number.}%
\end{array}%
\right.  \label{eq1}
\end{equation}%

The following sequence is obtained via the composite function (iteration):

$\Lambda=\{ n,T(n),T(T(n)),T(T(T(n))),\cdots\}=\{ n,T(n),T^2(n),T^3(n),\cdots\}$.

Consequently, the Collatz conjecture can be stated as follows:

\textbf{Collatz conjecture 1:}  For any natural number $n$, there is finite natural number $m$, the sequence $\Lambda$ always leads to the integer $1$, namely $T^m(n)=1$.

The series $\Lambda$ is an infinite sequence of ultimately period $^{[3,4]}$. So we give another statement of the Collatz conjecture as the following.

\textbf{Collatz conjecture 2:} The series $\Lambda$ is an infinite sequence of ultimately period, the preperiod $\eta(n)$ varies with the initial value $n$, but the ultimately period is always $\{1,4,2\}$.

\section{ An algebra and graph representation of natural numbers}

\subsection{The Composition of odd-number and even-number functions}

A natural number is considered even if it can be divided by 2; if not, it is considered odd. According to the Peano's Axiom, the smallest natural number is 1. The set $N=\{1,2,3,\cdots \}$ of natural numbers can be divided into odd and even sets; in this paper, we will use the usual definition of natural numbers.

$ \{ \textbf{natural number} \}=\{ \textbf{odd number} \} \bigcup \{ \textbf{even number} \} $.

In the set of natural numbers where 1 is the smallest odd number and 2 is the smallest even number, we can use the expression $n=2k-1$ to indicate that it is an odd, and the expression $n=2k$ to indicate that it is an even, where $k$ is any natural number.

We introduce two functions $O(x)=2x+1$ to express odd numbers greater than 1, and $E(x)=2x$ to express even numbers, where $x$ is any natural number in $N$.

We define a strictly increase monotonically piecewise function $f(n)$, from a natural number $n$ it generates two cases: odd or even numbers:
\begin{equation}
f(n)=\left\{
\begin{array}{cc}
2n+1=O(n), & \text{ the value is odd number,} \\
2n~~~~=E(n), & \text{ the value is even number.}%
\end{array}%
\right.  \label{eq1}
\end{equation}%

{\bf Definition 1}\quad A natural number $n$ is obtained by composition of \emph{the odd-number function} $O(x)=2x+1$ and \emph{the even-number function} $E(x)=2x$ several times, namely
\[
 n = f(f (\cdots f(1))) = f ^{k} (1),
\]
the function $f$ is either odd-number function $O(x)$ or even-number function $E(x)$.

For example, $f(1)=O(1)=3, f(1)=E(1)=2$, \\
$7 =f^2(1)= 2 \cdot 3 + 1 = 2 \cdot (2 \cdot 1 + 1) + 1) = O(O(1))$,\\
$189=2 \cdot 94+1= 2 \cdot (2 \cdot (47)+1 = 2 \cdot (2 \cdot (2 \cdot 23+1))+1 = 2 \cdot (2 \cdot (2 \cdot (2 \cdot 11+1)+1))+1  \\
= 2 \cdot (2 \cdot (2 \cdot (2 \cdot (2 \cdot 5+1)+1)+1))+1= 2 \cdot (2 \cdot (2 \cdot (2 \cdot (2 \cdot (2\cdot 2+1)+1)+1)+1))+1 \\
=2 \cdot (2 \cdot (2 \cdot (2 \cdot (2 \cdot (2\cdot (2\cdot 1)+1)+1)+1)+1))+1 = f^{7}(1)= O(E(O(O(O(O(E(1)))) $.

Any natural integer $n$ is the value of the finite times composite function of the odd-number and even-number functions beginning at 1. In particular, If $n=f(f(\cdots f(1)))=f^r(1)$, then $f^{-r}(n)=1$ is the inverse function.

\subsection{Use binary string to represent natural numbers}

Using binary string to represent a natural number $n$, we can more clearly express the odd-number and even-number functions starting at 1 composite process of a natural number. The string indicates the order of composition of $O(x)$ and $E(x)$; 0 implicates an even-number function, and 1 implicates an odd-number function.

A natural number's binary string represents its odd-even composite function; from left to right, the $i(i > 0)$ odd-number function $O(x)$ is represented by the $1$ in the $i(i > 0)$-bit, and the equivalent even-number function $E(x)$ is represented by $0$.

To convert a given natural number n to its binary string, we can use the following recursion steps:

1. Iterate over the number n, repeatedly dividing it by 2 and keeping track of the remainder.

2. Append the remainder to the binary string.

3. Continue dividing n by 2 until it reaches 0.

4. Reverse the binary string obtained from step 2 to get the correct binary representation of the number n.

For instance, the procedure of the composite function of 60 is displayed in Fig. 1.

\begin{figure*}[htbp] % 或者 figure 环境，取决于你的需求
  \centering
  \includegraphics[width=340pt]{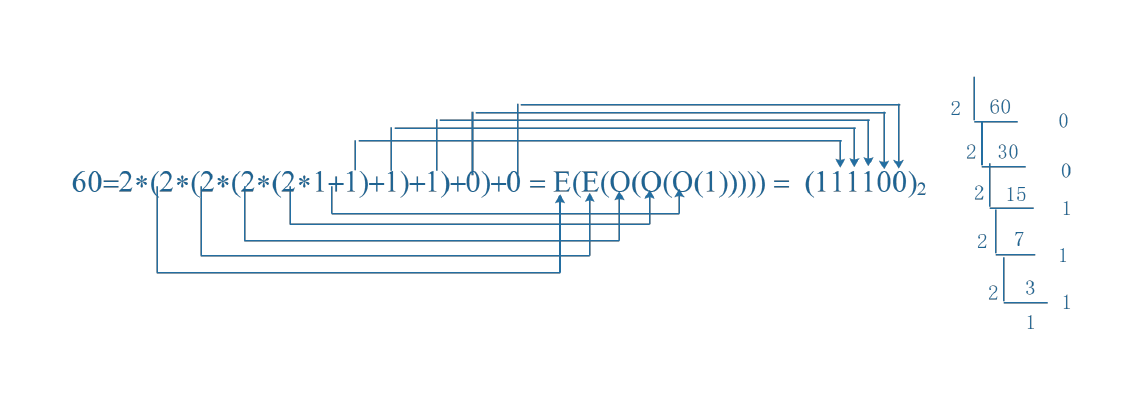} % 假设你使用 pt 作为单位
  %\captionsetup{font=small} % 如果需要设置标题字体大小为 small
  \caption{Natural number $60=(111100)_2$ is obtained starting 1 through
  the composition of five even-number and odd-number functions.}\label{fig:fig1}
\end{figure*}

\subsection{Use a graph to represent the composite procedure}
In order to give an intuitive impression, we provide \textbf{a full binary directed tree} to represent the procedure of composite function of the odd-functions $O(x)$ and (or) even-functions $E(x)$ of a given natural number $n$, the root is the smallest number 1. For per vertex, its left-child is an even number which double itself, its binary string is appended by 0, right-child is an odd number which double itself and add 1, its binary string is appended by 1. The full binary directed tree, as in Fig. 2, is a very good representation of some natural numbers.

\textbf{Proposition 1} A binary string's length indicates its level in the full binary directed tree, and a binary string's length minus one represents the number of times the odd- and even-number composite functions occur.

\begin{figure*} [htbp]
  \centering
 \includegraphics[width=360pt]{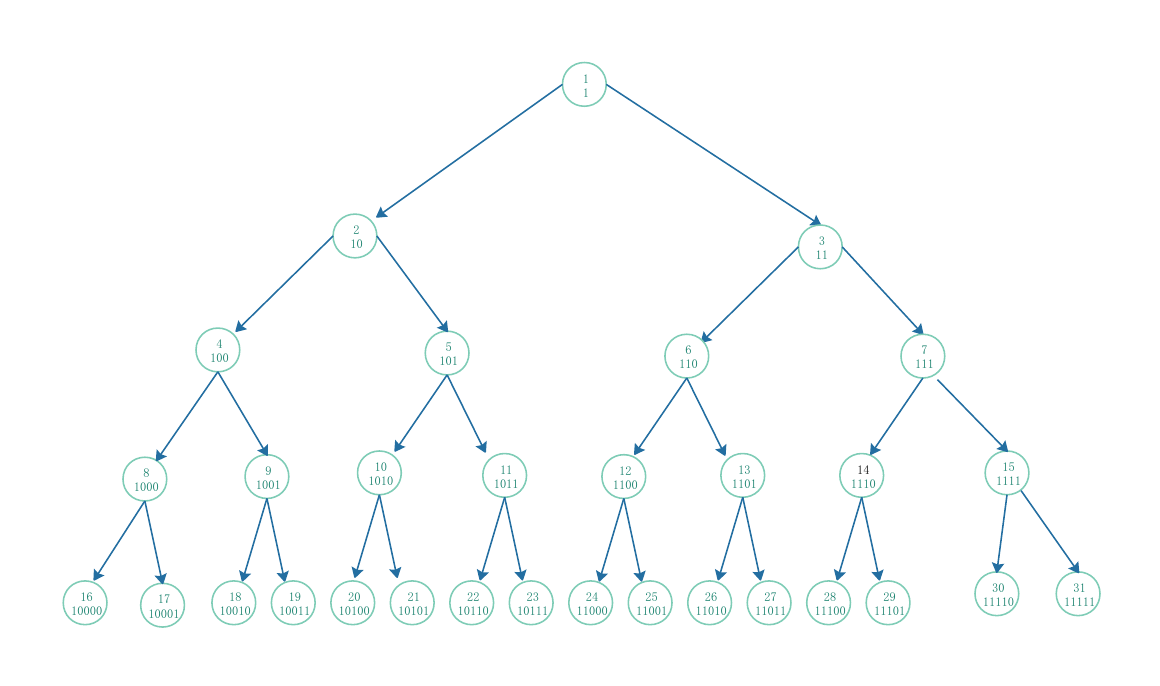}
  \centerline {\small Fig. 2 \quad The representation of natural number set is a full binary directed tree.}
  \label{fig:fig1}
\end{figure*}

For a natural number $n$, its binary string, composed of 0s and 1s from left to right, represents the path starting from the root node 1, and tracing down to the current node $n$ in the full binary directed tree. In this tree, each node only has one path leading from the root to itself. As an illustration, the procedures for the composite odd number 21 and the even number 28 are depicted in Figure 3.
\begin{figure*} [htbp]
  \centering
 \includegraphics[width=320pt]{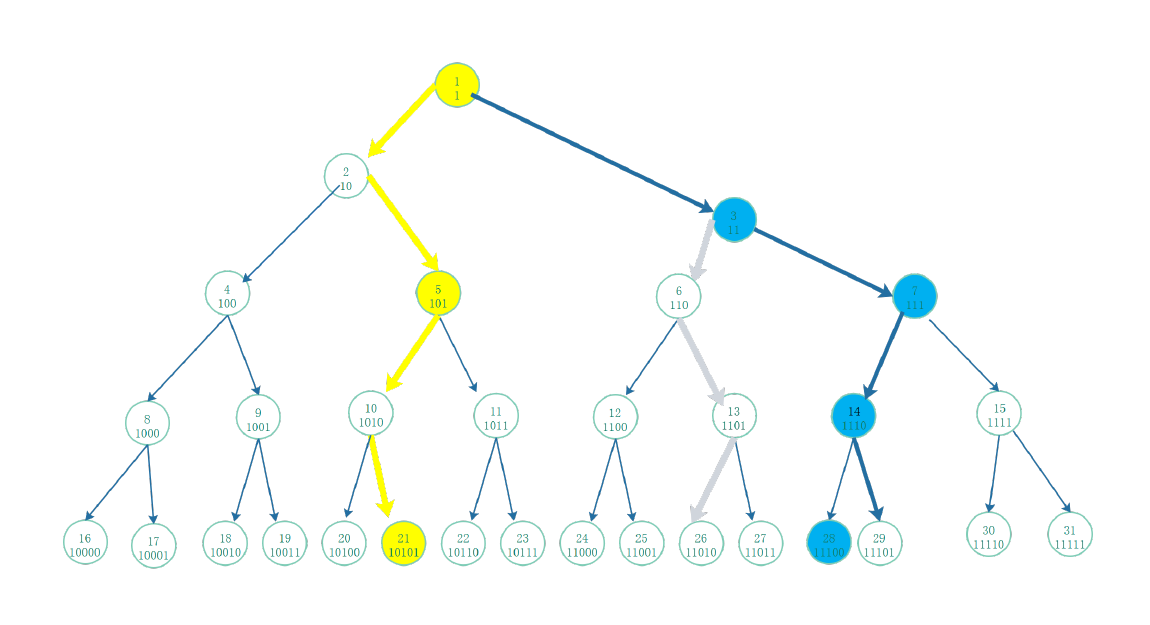}
  \centerline {\small Fig. 3\quad $21=(10101)_2$ and $28=(11100)_2$ comes from the path from root 1 walk to}
  \centerline {\small  \quad   $10101$ and $11100$ accordingly appending 1 or 0 to the nodes in succession.}
  \label{fig:fig1}
\end{figure*}

\subsection{Another partition of the natural number set}

To analyze the iterative functions of the odd-function and even-function, we introduce a new partition of the natural number set within the full binary directed tree. We designate distinct names to the numbers located on the left path, right path, and between the two paths in the tree.

We give the definitions of three kinds of natural numbers:

{\bf Definition 2}{\bf(i)} A natural number, $ O^m(1)=2^m-1=2^{m-1}+2^{m-2}+\cdots+2+1=(11\cdots 1)_2$, is obtained by applying the odd-number function $O(x)$ $m$ components. We call it as \textbf{pure odd number}. For instance, those are pure odd numbers: $ 3=(11)_2, 7=(111)_2, 15=(1111)_2, 31=(11111)_2, 63=(111111)_2,\cdots $. These are located in the full binary directed tree of Figure 2, which is the right leg of the isosceles triangle.

{\bf(ii)} A natural number, $ E^m(1)=2^m=(10\cdots 0)_2 $ , is obtained by applying the even-number function $E(x)$ $m$ components. We call it as \textbf{pure even number}. For instance, those are pure even numbers: $ 2=(10)_2 $, $4=(100)_2 $,  $8=(1000)_2 $,  $16=(10000)_2 $, $32=(100000)_2 $,  $64=(1000000)_2 $. Those are located in the left leg of the isosceles triangle, namely the full binary directed tree of the Figure 2.

{\bf(iii)} The natural number obtained by the composition of odd function $O(x)$ and even function $E(x)$, we call it \textbf{mixed number}. Such as, $ 18=(10010)_2, 28=(11100)_2, 67=(1000011)_2, 309=(100110101)_2 $. Those are in the inside of the isosceles triangle, the full binary directed tree of the Figure 2.

In particular, the natural numbers obtained by the finite alternately composition of the odd function $O(x)$ and the even function $E(x)$, namely, $[E(O(1))]^m=(101\cdots 101)_2$. Such as
$\ 5=(101)_2, 21=(10101)_2, 85=(1010101)_2, 341=(101010101)_2,\\
\  1365=(10101010101)_2, 5461=(1010101010101)_2,\cdots . \\ $ We call \textbf{hard numbers}.

The traversal path in the full binary directed tree from the root down along the arcs, for each natural integer $n$, is its binary string $1\times \times$, where the left-child appended $0$ for each node is an even number and the right-child appended $1$ for each node is an odd number. For instance, in Fig. 3, $21=(10101)_2$ originates at the root 1 and proceeds down 2,5,10, ultimately reaching 21. To the nodes, $1\rightarrow 10\rightarrow 101\rightarrow 1010 \rightarrow 10101$, 0,1,0,1 are appended. In addition, for $28=(11100)_2$, the appendix 1,1,0,0 is added to the nodes, $1\rightarrow 11\rightarrow 111\rightarrow 1110 \rightarrow 11100$, accordingly in dicimal it traces from the root 1 down 3,7,14, and ultimately reaches 28.

{\bf Property 2}\quad The set of natural numbers can be partited into three different sets:

\textbf{\{natural number\}}=

\ \ \ \ \ \ \ \ \ \ \ \ \ \ \ \ \ \{\textbf{pure even number}\}$\cup $\{%
\textbf{pure odd number}\}$\cup $\{\textbf{mixed number}\},
where
\{\textbf{mixed number\}=}\{\textbf{mixed even number}\}$\cup $\{\textbf{mixed odd number}\}

{\bf Example 1} (1) $60$, $97$ are mixed numbers.

(2)$64,1180591620717411303424$ are pure even numbers.

(3)$63,1180591620717411303423$ are pure odd numbers.

When we convert those natural numbers from decimal to binary, the facts are obvious.

(1) $60=(111100)_2 $ is a mixed-even number, $97=(1100001)_2$ is a mixed-odd number.

(2) $64=2^6=(1000000)$, $1180591620717411303424=2^{70}=(10000\dots 0)$ are pure even numbers:

$63=(111111)_2$, $1180591620717411303423=2^{70}-1=(11\dots 1)$ are pure odd numbers.

\section{Two functions compare} % (fold)

\subsection{The piecewise and iterative functions}

A piecewise function is a mathematical function that is defined by different rules or formulas over different intervals or regions of its domain. Piecewise functions are commonly used to model situations where different rules apply in different circumstances or to account for discontinuities in a function's behavior. The Collatz function, denoted by $T(n)$, can be expressed as a piecewise function, with separate cases for odd and even numbers.

An iterative function is a function that is repeatedly applied to its own output. In other words, the output of the function is used as the input for the next iteration of the function. Iterative methods involve using iterate functions to repeatedly update an initial estimate or solution until a desired level of accuracy is achieved. Iteration means repetition, and with more repetition, things will change in nature.

In order to proof the Collatz conjecture 1 in section 1, finding the beingness and finiteness of the number $m$ in the expression $T^m(n)=1$ for a natural number $n$ is the main challenge. Iteration is the key to Collatz conjecture, and although there are only two cases where piecewise functions are combined with iterative functions, the result is difficult to control.

To visually represent the process of the iterative functions of odd- and even-number functions, i.e., $n=f^k(1)$, we propose a full binary directed tree. The following is the inverse functions $f^{-1}(x)$,
\begin{equation}
f^{-1}(n)=\left\{
\begin{array}{cc}
\frac{n-1}{2}, & \text{ If $n$ is odd number,} \\
\frac{n}{2}, & \text{ If $n$ is even number.}%
\end{array}%
\right.  \label{eq1}
\end{equation}%

For given natural $k$, the iterate formula $f^k(1)=n$ for a given natural number $n$, we know that $k$ is the length of the binary string of $n$ minus 1, it is also the level of the full binary directed tree. In decimal notation, we represent $n$, which obscures the composite process of odd- and even-number functions. When $n$ is represented as a binary string, it can be used to understand how odd- and even-number functions composite work.
\subsection{The Collatz function and reduced Collatz functions}

The iteration of the Collatz function is the key topic in discuss the proof procedure, we have the reduced Collatz function $^{[2,5,6]}$ $RT(x)$
\begin{equation}
RT(n)=\left\{
\begin{array}{cc}
\frac{3n+1}{2^m}=T^{m+1}(n), & \text{if }n\text{ is odd number, $\frac{3n+1}{2^m}$ is an odd number,} \\
\frac{n}{2^r}=T^r(n), & \text{if }n\text{ is even number, $\frac{n}{2^r}$ is an odd number.}%
\end{array}%
\right.  \label{eq3}
\end{equation}%

There are many different points for piecewise functions when comparing the Collatz function $T(x)$ with the function $f^{-1}(x)$, and the reduced Collatz function $RT(x)$ with the iteration function $f^k(x)=O(E(\dots (E(x)))$.

1) For any natural number $x$ the function $f^{-1}(x)$ and $f^{-r}(x)=(f^{-1}(n))^r$ are strictly monotonically decreasing.

2) The function $T(x)$ is rising in the case $x$ is an odd number, in the other case, is decrease.

3) The function $RT(x)$ describing the procedure of the iterative function of $T(x)$, is wavy when $x$ is a pure or mixed odd number, and decrease when $x$ is pure even or mixed even.

\begin{tabular}{|l|l|l|l|l|l|l|l|l|}
\hline
$RT$ & $\frac{3n+1}{2}$ & $\frac{3n+1}{2^{2}}$ & $\frac{3n+1}{2^{3}}$ & $%
\cdots $ & $\frac{n}{2}$ & $\frac{n}{2^{2}}$ & $\frac{n}{2^{3}}$ & $\cdots $
\\ \hline
$T$ & $T^{2}(n)$ & $T^{3}(n)$ & $T^{4}(n)$ & $\cdots $ & $T^{2}(n)$ & $%
T^{3}(n)$ & $T^{4}(n)$ & $\cdots $ \\ \hline
$monotonicity$ & $\uparrow \downarrow $ & $\uparrow \downarrow \downarrow $
& $\uparrow \downarrow \downarrow \downarrow $ & $\cdots $ & $\downarrow $ &
$\downarrow \downarrow \downarrow $ & $\downarrow \downarrow \downarrow
\downarrow $ & $\cdots $ \\ \hline
\end{tabular}

The function $RT(x)$ describes the procedure of iterative of $T(x)$. The wavy function is increasing at first, then goes through one or more decreasing processes, either as "increase -- decrease -- increase" or "increase -- decrease $\cdots$ decrease -- increase." For instance, the iterated sequence of Collatz functions is plotted in Figs. 4 and 5, where the beginning values are pure odd $255=2^8-1=(11111111)_2$ and mixed odd number $97=(1100001)_2$, respectively.
\begin{figure*} [htbp]
  \centering
  \includegraphics[width=200pt]{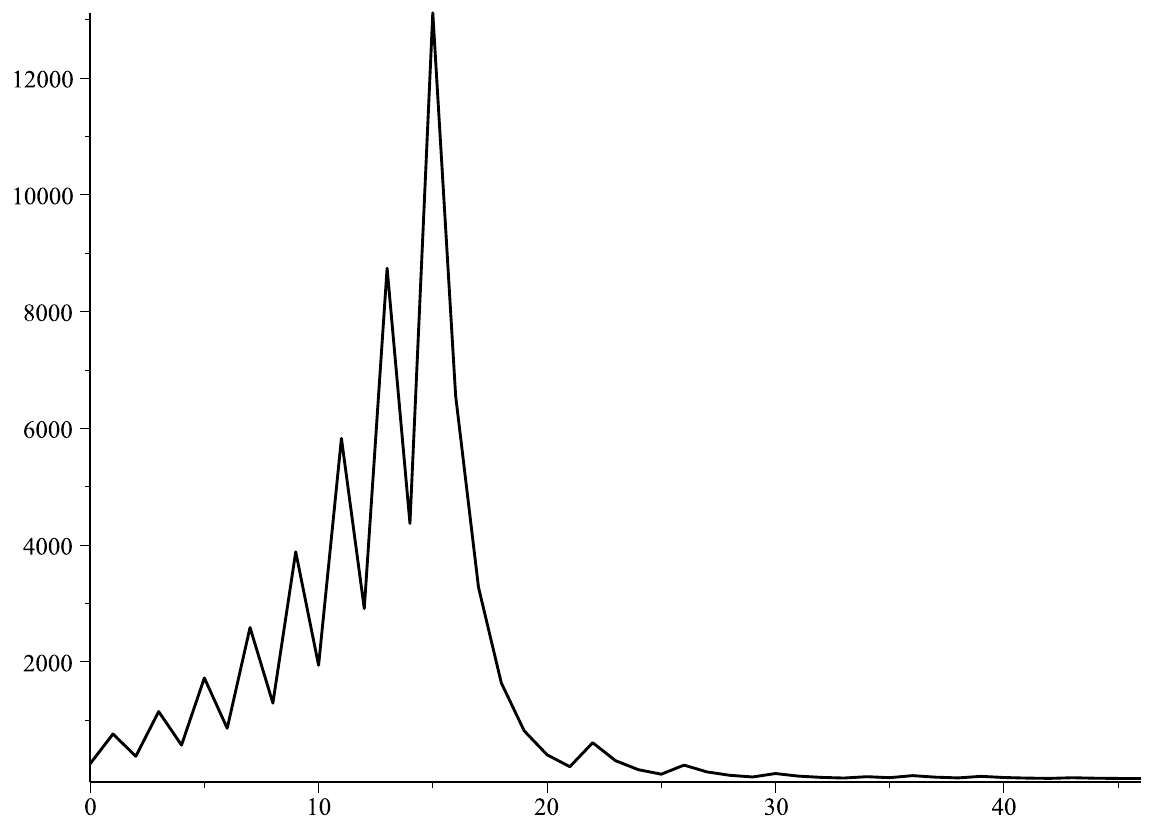}
  \centerline {\small Fig. 4\quad Point plot of a sequence of 47 iterations of the Collatz function for pure odd 255.}\label{fig:fig1}
\end{figure*}

\begin{figure*} [htbp]
  \centering
  \includegraphics[width=200pt]{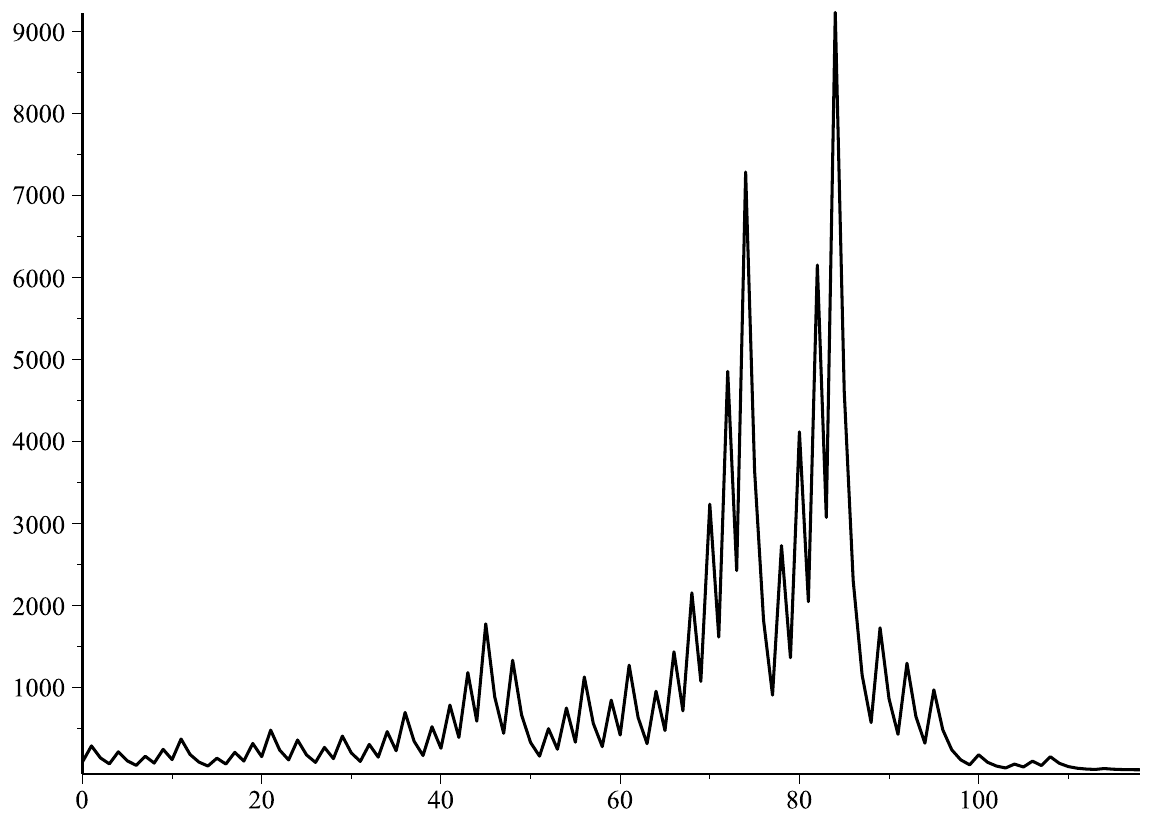}
  \centerline {\small Fig. 5\quad Point plot of a sequence of 118 iterations of the Collatz function for mixed odd 97.}\label{fig:fig1}
\end{figure*}

For a given natural number $n$, the Collatz iterative function converges to 1 in a finite number of steps and cycles indefinitely between the numbers 1, 4, and 2 for an infinite number of iterations.

\section{Using binary string to explore the Collatz conjecture} % (fold)
\label{sec:Proof the Collatz conjecture}

If the Collatz function (1) is expressed in binary form as
\begin{equation}
T(n)=\left\{
\begin{array}{cc}
(11)_{2}\cdot (1\times \cdots \times 1)_{2}+1=(1\times \times \times
10\cdots 0)_{2}, & \text{if }n\text{ is odd number,} \\
\frac{(1\times \cdots \times 10)_{2}}{(10)_{2}}=(1\times \cdots
\times 1)_{2}, & \text{if }n\text{ is even number.}%
\end{array}%
\right.  \label{a2}
\end{equation}
The characteristics of the left side and right side and the penultimate of the binary string are illustrated by the fig. 6.
\begin{figure*} [htbp]
  \centering
  \includegraphics[width=340pt]{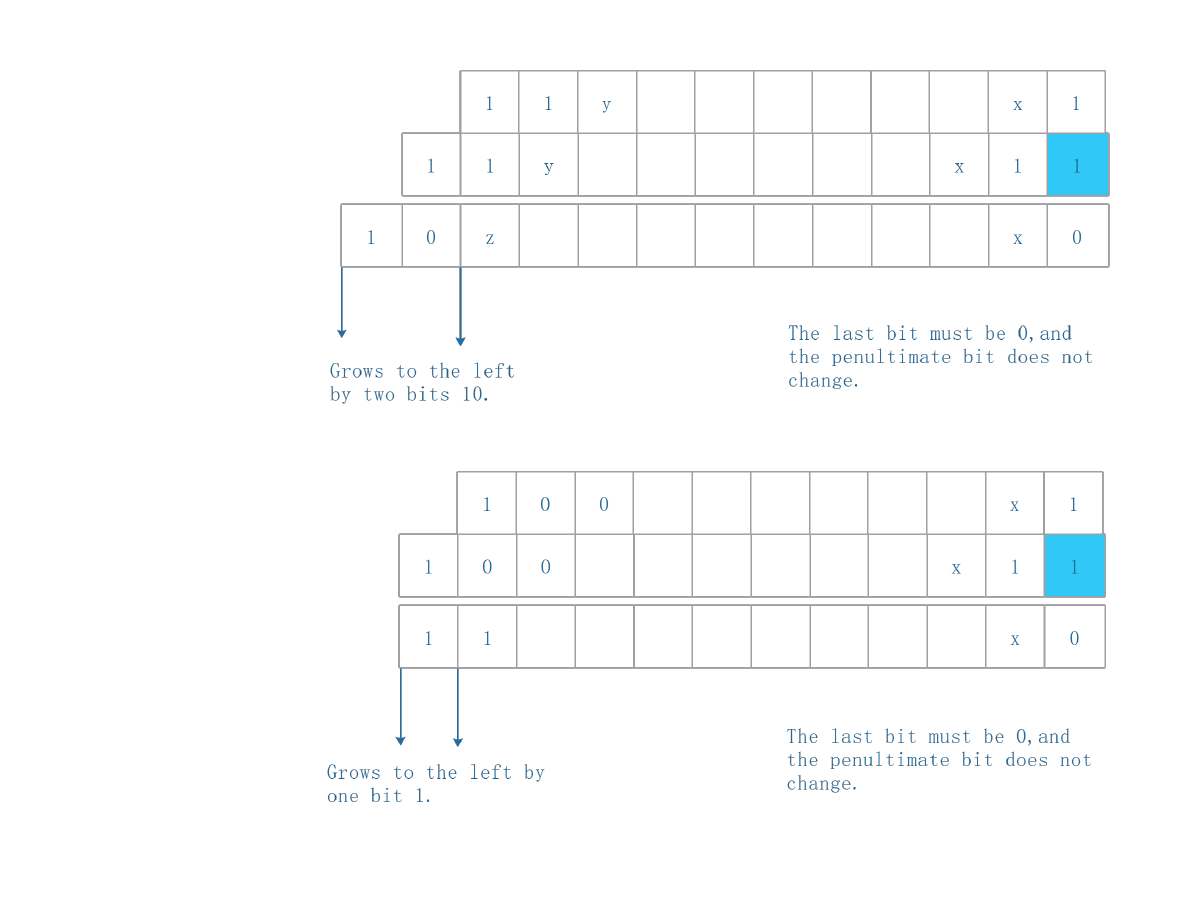}
  \centerline {\small Fig. 6\quad Binary representation of the Collatz function,}
  \centerline{ grows by appending 1 or 10 to the left side of the binary string}
  \centerline{ delete at least one 0s in the right side of the binary string.} \label{fig:fig1}
\end{figure*}

Then we use binary string to illustrate the reduced Collatz function (4) as the follows,

\begin{equation}
RT(n)=\left\{
\begin{array}{cc}
\frac{(1\times \times \times 10\cdots 0)_{2}}{(10\cdots 0)_{2}}=(1\times \times \times 1)_{2}, & \text{if }n\text{ is an odd, $\frac{3n+1}{2^m}$ is an odd,} \\
\frac{(1\times \times \times 10\cdots 0)_{2}}{(10\cdots 0)_{2}}=(1\times \times \times 1)_{2} & \text{if }n\text{ is an even, $\frac{n}{2^r}$ is an odd.}%
\end{array}%
\right.  \label{eq3}
\end{equation}%
where $\times $ is 0 or 1.

We adopt the binary representation method for specific natural numbers and use mathematical experimental methods to obtain their Collatz sequences. The following are three kinds forms to describe the Collatz sequences respectively: (A)tabular, (B)scratch paper and (C)mathematical experiment, for.
We use the binary representation method for specific natural numbers and employ mathematical experimental techniques to generate their Collatz sequences. The Collatz sequences for the numbers $n=10027, 255, 78736985$ are presented in three different formats: (A) tabular, (B) scratch paper, and (C) mathematical experimentation.

(A) For the formula $T^{91}(10027)=1$, we apply the mathematical software Maple get the sequence $T^i(10027), i=0..210$ are in decimal and binary as the follows.
\begin{tabular}{|r|r|}
\hline
10027=(10011100101011)$_{2}\rightarrow $ & (111010110000010)$_{2}\rightarrow
$ \\ \hline
15041=(11101011000001)$_{2}\rightarrow $ & (1011000001000100)$%
_{2}\rightarrow $ \\ \hline
11281=(10110000010001)$_{2}\rightarrow $ & (1000010000110100)$%
_{2}\rightarrow $ \\ \hline
8461=(10000100001101)$_{2}\rightarrow $ & (110001100101000)$_{2}\rightarrow $
\\ \hline
3173=(110001100101)$_{2}\rightarrow $ & (10010100110000)$_{2}\rightarrow $
\\ \hline
595=(1001010011)$_{2}\rightarrow $ & (11011111010)$_{2}\rightarrow $ \\
\hline
893=(1101111101)$_{2}\rightarrow $ & (101001111000)$_{2}\rightarrow $ \\
\hline
335=(101001111)$_{2}\rightarrow $ & (1111101110)$_{2}\rightarrow $ \\ \hline
503=(111110111)$_{2}\rightarrow $ & (10111100110)$_{2}\rightarrow $ \\ \hline
755=(1011110011)$_{2}\rightarrow $ & (100011011010)$_{2}\rightarrow $ \\
\hline
1133=(10001101101)$_{2}\rightarrow $ & (110101001000)$_{2}\rightarrow $ \\
\hline
425=(110101001)$_{2}\rightarrow $ & (10011111100)$_{2}\rightarrow $ \\ \hline
319=(100111111)$_{2}\rightarrow $ & (1110111110)$_{2}\rightarrow $ \\ \hline
479=(111011111)$_{2}\rightarrow $ & (10110011110)$_{2}\rightarrow $ \\ \hline
719=(1011001111)$_{2}\rightarrow $ & (100001101110)$_{2}\rightarrow $ \\
\hline
1079=(10000110111)$_{2}\rightarrow $ & (110010100110)$_{2}\rightarrow $ \\
\hline
1619=(11001010011)$_{2}\rightarrow $ & (1001011111010)$_{2}\rightarrow $ \\
\hline
2429=(100101111101)$_{2}\rightarrow $ & (1110001111000)$_{2}\rightarrow $ \\
\hline
911=(1110001111)$_{2}\rightarrow $ & (101010101110)$_{2}\rightarrow $ \\
\hline
1367=(10101010111)$_{2}\rightarrow $ & (1000000000110)$_{2}\rightarrow $ \\
\hline
2051=(100000000011)$_{2}\rightarrow $ & (1100000001010)$_{2}\rightarrow $ \\
\hline
3077=(110000000101)$_{2}\rightarrow $ & (10010000010000)$_{2}\rightarrow $
\\ \hline
577=(1001000001)$_{2}\rightarrow $ & (11011000100)$_{2}\rightarrow $ \\
\hline
433=(110110001)$_{2}\rightarrow $ & (10100010100)$_{2}\rightarrow $ \\ \hline
325=(101000101)$_{2}\rightarrow $ & (1111010000)$_{2}\rightarrow $ \\ \hline
61=(111101)$_{2}\rightarrow $ & (10111000)$_{2}\rightarrow $ \\ \hline
23=(10111)$_{2}\rightarrow $ & (1000110)$_{2}\rightarrow $ \\ \hline
35=(100011)$_{2}\rightarrow $ & (1101010)$_{2}\rightarrow $ \\ \hline
53=(110101)$_{2}\rightarrow $ & (10100000)$_{2}\rightarrow $ \\ \hline
5=(101)$_{2}\rightarrow $ & (10000)$_{2}\rightarrow $ \\ \hline
\end{tabular}

(B) For the formula $T^{47}(255)=1$, we get the sequence $T^i(255), i=0..47$ are in decimal and binary as the following scratch paper as in Fig. 7. The diagonal arrow line means that the length of the trailing field is successively reduced by one digit until it is 1, the red vertical arrow line means that the pure even number is successively divided by 2 until it is odd, and the red horizontal arrow line means the 3n+1 operation.

\begin{figure*} [htbp]
  \centering
  \includegraphics[width=340pt]{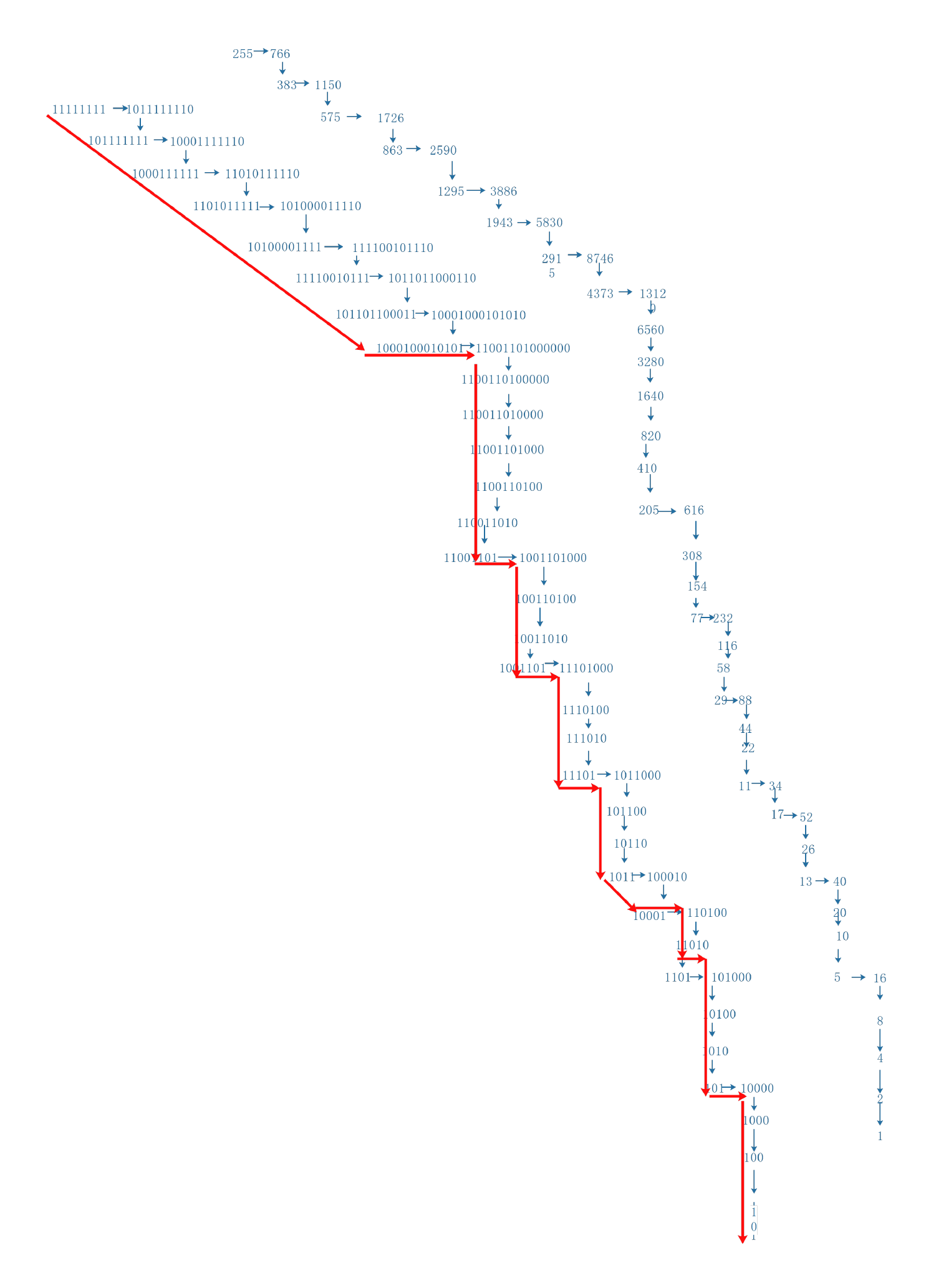}
  \centerline {\small Fig. 7\quad The scratch paper of a sequence of 47 iterations of the}
  \centerline {\small Collatz function for pure odd 255 in decimal and binary forms.}\label{fig:fig1}
\end{figure*}

(C) For the procedure of $T^{210}(78736985)=1$ from 0 to 210, we can see Appendix to this article for details.

We observe the procedure of iterative Collatz function, namely the reduced Collatz function (6), i.e., the Collatz sequences, that we pay close attention to the zeros in right side of binary strings of an even number and the \emph{end-substring} between the first 0 encountered from right to left which is make of 1. For instance, for $1011001$ the end-substring is $1$,  for $1011001111$ the end-substring is $1111$,  for $11111$ the end-substring is itself $11111$. There are many properties of the end-substrings.

In the Collatz sequence represented by a binary string, looking backwards from $n$:

(1) If there are several zeros at the end, remove one at a time until all zeros are deleted, and the number becomes odd.

(2) When the number of bits at the end-substring is more than 1, the adjacent binary string must have only a single zero at the end. Delete this zero to make it the next odd number, and continue these steps until only a 1 is left at the end.

(3) If the number of digits at the end-substring is only one, the adjacent binary string must end with several zeros. Deleting these zeros in sequence will result in the following two scenarios at the end of the binary string: (i) One bit 1, (ii) More than one bit 1.

\section{Proof the Collatz conjecture} % (fold)
\label{sec:Proof the Collatz conjecture}

We have known that mathematics formula about geometric progression with initial term 1 and common ratio $x$, the sum of the first $k$ terms is
\begin{equation}
x^k+x^{k-1}+\cdots+x+1=\frac{x^{k+1}-1}{x-1}
\end{equation}%

when $x=2$, there are two formulas

\begin{equation}
2^{k-1}+2^{k-2}+\cdots 2+1=2^{k}-1
\end{equation}%
\begin{equation}
2^{k}+2^{k}=2^{k+1}%
\end{equation}%

The substantive characteristics is that the powers of $2$ must be continuous natural numbers, this is the key to our proof method to solve the Collatz conjecture.

\begin{proof} For a given natural number $n$, it can always be represented by a binary string, that is, it must be found on the full binary directed tree we mentioned earlier, which means that its structure can be determined by the difference between the corresponding powers in the binary string representation. Every time the Collatz function operation is performed, its structure is adjusted accordingly, and $3n+1$ can add a term $2^{r}$ between the corresponding two terms $2^{r+1}, 2^{m+1}$. Similarly, $\frac{3n+1}{2^h}$ can change the power of each bit term by term. This ensures that the bits in its binary string will definitely change.

(i) For a given natural number $n=2^k$ is a pure even, then the smallest natural number 1 can be reached by simply repeating the $k$ times Collatz function divided by 2.

If $n=(1\times 10\cdots0)_2=2^m+\cdots+2^k$ is a mixed even number, one can delete all zeros in end, it becomes as a mixed odd $\frac{n}{2^k}=(1\times 1)_2=2^{m-k}+\cdots+2^0$.

(ii) when a given natural number $n$ is either pure odd or mixed odd, its binary string as $ n=(1\times\cdots\times 1)_2=2^r+2^m+\cdots +1 $, then
\begin{flalign*}
  & 3n+1=2n+(n+1)  \\
  &=2^{r+1}+2^{m+1}+\cdots+2+2^r+2^m+\cdots+1+1 \\
 &    =2^{r+1}+2^{m+1}+\cdots+2^0+2^r+2^m+\cdots+2^0  \\
 &=2^{r+1}+2^r+2^{m+1}+2^m+\cdots+\cdots+2^h
\end{flalign*}

 If the length of the end-substring of $n$ is 1, the length of end-substring of the binary string $3x+1$ is either 1 or bigger than 1.

Two formulas (8) and (9) can be used to modify the structure of the binary string $n$ to the binary string $3n+1$. That is, there is an appended term $2^{r-h}$ in two equivalent terms, $2^{r+1-h}$ and $2^{m+1-h}$. When the zeros in the middle of a binary string are compared to a bubble, it means that these zeros are gradually being driven out of the rightmost end by $3n+1$. It is the same as progressively removing the bubbles hidden in the sponge using a means $3n+1$. Once 3n+1, 0 shifts one bit to the right, i.e., the length of the associated binary substring is reduced by one bit, when the length of the end-substring is greater than 1. The end-substring length of the binary string $3x+1$ is either greater than or equal to 1 if the length of the end-substring of $n$ is 1.

We shall then divide $3n+1$ by the last term $2^h$,
\begin{flalign*}
 &  \frac{3n+1}{2^h}  =2^{r+1-h}+2^{r-h}+2^{m+1-h}+2^{m-h}+\cdots+\cdots+2^0 \\
\end{flalign*}
get another odd number, this is the value of reduced Collatz function (4) or (6). And so on, finitely steps after finally we get a pure even number $2^t$, this is the case in above (i), thus the Collatz conjecture hold on.

We illustrate the procedure by a mixed odd number $n=67$ and hard number set in the following,
\begin{align*}
&67=(1000011)_2=2^6+2^1+2^0  \\
&3\cdot67+1=2\cdot67+67+1 =2^7+2^2+2^1+2^6+2^1+2^0+2^0=2^7+2^6+2^3+2^1 \\
&3\cdot101+1=2\cdot101+101+1=2^7+2^6+2^3+2^1+2^6+2^5+2^2+2^0+2^0=2^8+2^5+2^4 \\
&3\cdot19+1=2\cdot19+19+1=2^5+2^2+2^1+2^4+2^1+2^0+2^0=2^5+2^4+2^3+2^1 \\
&3\cdot29+1=2\cdot29+29+1=2^5+2^4+2^3+2^1+2^4+2^3+2^2+2^0=2^6+2^4+2^3 \\
&3\cdot11+1=2\cdot11+11+1=2^4+2^2+2^1+2^3+2^1+2^0+2^0=2^5+2^1 \\
&3\cdot17+1=2\cdot17+17+1=2^5+2^1+2^4+2^0+2^0=2^5+2^4+2^2 \\
&3\cdot13+1=2\cdot13+13+1=2^4+2^3+2^1+2^3+2^2+2^0+2^0=2^5+2^3 \\
&3\cdot5+1=2\cdot5+5+1=2^3+2^1+2^2+2^0+2^0=2^4 \\
&1
\end{align*}

For a special class of mixed numbers, the hard number $\frac{4^k-1}{3}=(101\cdots101)_2$, then its Collatz sequent result is
\begin{align*}
&a_{k}=\frac{4^{k}-1}{3}=\frac{4^{k}-1}{4-1}=4^{k-1}+4^{k-2}+\cdots+4+1=(101\cdots 101\cdots 101)_{2},\\
&T(a_{k})=3a_{k}+1=4^k=2^{2k}=(10\cdots 0)_{2},T^{2k+1}(a_{k})=1.
\end{align*}%

This means that the Collatz conjecture is valid for this case. Therefore we have proved the Collatz conjecture 1 at section 1 of this paper.\end{proof}

\begin{proof}
For the smallest number $1=2^0$
\begin{align*}
 &T(1)=3\times 1+1=2\times 1+1+1=2^1+2^0+2^0=2^1+2^1=2^2=4,\\
 &T^2(1)=2,\\
 &T^3(1)=1,
\end{align*}
 and so on, this is a cycle $\{1,4,2\}$, This is the Collatz conjecture 2 at section 1 of this paper.\end{proof}

\section{Conclusion}

From previous proof of the conjecture, it becomes a theorem.

\textbf{Theorem 1} There exists a finite natural number $m$ for every natural number $n$. The Collatz sequence always arrives at the integer $1$, that is, $T^m(n)=1$.

\textbf{Theorem 2} For any positive integer $n$, the sequence of the iteration of the Collatz function is an ultimately periodic sequence, its preperiod $\eta(n)$ is a related-to $n$, and the least period is $\{1,4,2\}$.

\section{Appendix}
\includepdf[pages=-]{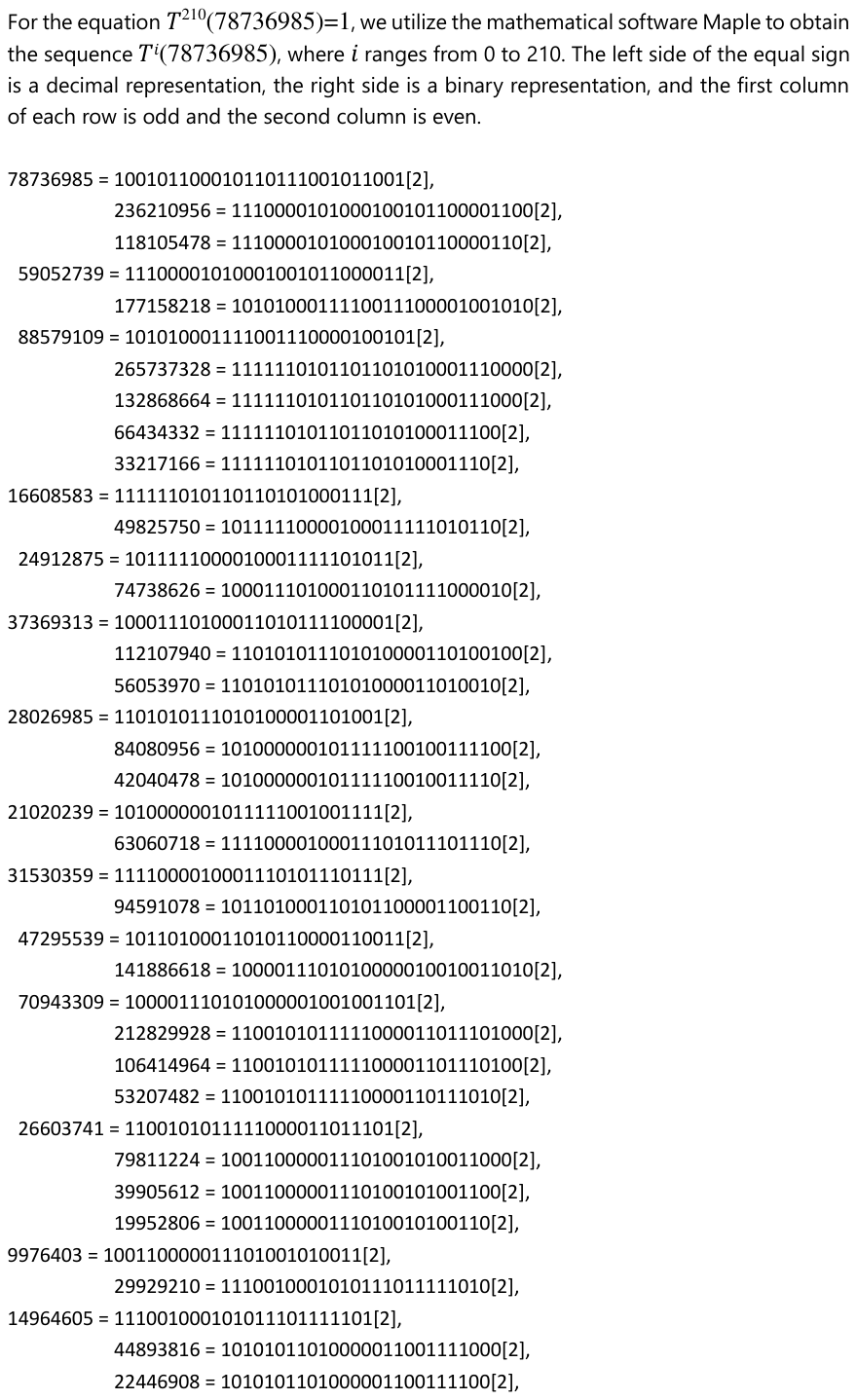}

\end{document}